\documentclass{article}
\usepackage{verbatim}
\usepackage[utf8]{inputenc}
\usepackage{amsmath,amssymb,amsthm}
\usepackage{tikz-cd}
\usepackage{bbm}

\usepackage{stmaryrd}
\usepackage[new]{old-arrows}
\usepackage[colorinlistoftodos]{todonotes}
\usepackage{graphicx}
\usepackage{hyperref}

\newcommand{\kk}{ \mathbbm{k} }
\newcommand{\AAA}{\mathbb A}

\newcommand{\lbm}{\left[ \begin{matrix}}
\newcommand{\rem}{\end{matrix} \right]}
\newcommand{\SL}{\sum\limits}

\newcommand{\ra}{\rightarrow}

\newcommand{\lp}{\left(}
\newcommand{\rp}{\right)}

\newcommand{\lav}{\left|}
\newcommand{\rav}{\right|}

\newcommand{\np}{\newpage  \noindent}

\newcommand{\p}{\partial}

\theoremstyle{definition}

\begin{document}

\title{"About the Solution of the General Equations of Fifth and Sixth Degree (Excerpt from a letter to Mr. K. Hensel.)"\footnote{Printed from the Dirichletbande (Bd. 129) of the Journal for Pure and Applied Mathematics.}\\
\hrulefill \\
English Translation of "\"Uber die Aufl\"osung der allgemeinen Gleichungen f\"unften und sechsten Grades (Auszug aus einem Schreiben an Herrn K. Hensel)" \footnote{Abgedruckt aus dem Dirichletbande (Bd. 129) des Journals für reine und angewandte Matiiematik.}}
\author{Original Author: Felix Klein\\ Translation: Alex Sutherland}
\date{Original Year of Publication: 1905\\ Year of Translation: 2019}

\maketitle
\begin{center}
    {\Large First appeared in the \textit{Dirichletbande} (Vol. 129) of the Journal for Pure and Applied Mathematics (1905). Reprinted in the Math. Annalen, Vol. 61 (1905).}
\end{center}

\vfill

\noindent
There are a few footnotes added to this translation, that may add extra insight or clarification to the reader. All of these footnotes begin with "Translator's Note:". 

\begin{center}
    \href{https://bit.ly/2HK6gdM}{Here is a link to the original paper.}
\end{center}

\newpage

\section{Introduction}
From,\\
\\
Felix Klein in G\"ottingen.\\

By responding to your earnest request to contribute to the journal's book dedicated to the memory of Dirichlet, I refer to a note I published six years ago in the \textit{Rendiconti dell'Accademia dei Lincei} \cite{Klein1899}, and in which I outlined a general solution of equations of sixth degree. \\

I set myself the goal of explaining in more detail and in more concrete terms what was suggested there. In fact, even an expert of the relevant literature (such as Mr. Lachtin) has not taken the approach in question in its simplicity (as I will explain more below). \footnote{1901, Moscow Mathematical Collection, Vol. XXLI, pp. 181-218 (Russian) \\} \\

Moreover, I act under the impulses of my old friend Mr. Gordan, who has recently turned his great algebraic ability to the problem in question. Mr. Gordan will soon publish a first relevant treatise in the Math. Annalen \cite{Gordan} \footnote{A contribution to the solution of the general equations of the sixth degree. Compare with a message to the Heidelberg International Congress of Mathematics \\}. But this is only a beginning; I hope that his continued efforts will succeed in clarifying the subject in every respect as fully as we have been able to do in the past with the theory of equations of the fifth degree.\\

I would like to discuss this theory of the equations of the fifth degree in advance, as I summarized them in my "Lectures on the Icosahedron" \cite{KleinLectures}, in such a way that I bring forth those moments which are generalized when dealing with the considerations of the sixth degree. In Chapter V of these Lectures, I have dealt with two methods for solving equations of the fifth degree (which, incidentally, differ only by the order in which the steps are carried out) and the second of these methods will prove to be the natural continuation of Kronecker's (and Brioschi's) method. This method, like the first one, is developed in geometric form, with special relations that emerge only in equations of the fifth degree. Instead, I refer here to the algebraic justification developed in Volume 15 of the Math. Annalen \cite{Klein1879} and accompanied with reflections on the solution of arbitrary higher equations. \footnote{In particular, see Section 4 - "The formulas of Kronecker and Brioschi for the fifth degree"} \\

The icosahedral theory of the equations of the fifth degree and the general considerations connected with it have since been portrayed several times by others, especially in the second volume of the excellent textbook on algebra by Mr. Weber \cite{Weber}, as well as in the detailed report that Mr. Wiman has made in Volume I of the Encyclopedia of Mathematical Sciences \cite{Wiman1898}. Nevertheless, it seems that the basic meaning of the whole approach in the mathematical audience is still often not understood. It is not a matter of considerations which are to the sides of the earlier investigations on the solution of equations of the fifth degree, but of those which claim to constitute the very core of these earlier investigations. Accordingly, in the following report, I will try to describe the main points of the theory (which will later be found \textit{mutatis mutandis} in the approach to the equations of the sixth degree) as accurately as possible while maintaining brevity. \\

The first is that we have the icosahedral equation, i.e. the equation of the sixtieth degree, which is written in the above Lectures as follows:
\begin{align}
    \frac{H^3(x)}{1728f^5(x)} = X
\end{align}

\noindent
as a \textit{Normalgleichung sui generis} (Normal General Equation) which, by virtue of their excellent qualities, is the next generalization of the "pure" equations:
\begin{align}
    x^n = X.
\end{align}

\noindent
In fact, given any root of (1), the 60 roots of (1) can be calculated by the 60 linear equations that are already known (the icosahedral substitutions), just as the $n$ roots of (2) can be found from any one of them by the $n$ substitutions given by $x' = e^\frac{2\pi ik}{n}x$. Now, the group of icosahedral substitutions proves to be isomorphic with the group of 60 alternating permutations on five letters (i.e. $A_5$). In this way, it is impossible to trace the solution of the general equations of the fifth degree back to a sequence of pure equations (2). The task is thus \textit{to solve the equations of the fifth degree with the help of an icosahedral equation.} Here we distinguish an algebraic and a transcendental part of the investigation. The first part will deal with the algebraic construction of a root $x$ of an icosahedral equation (1) from the roots $z_1,\dotsc,z_5$ of a given fifth degree equation -  the parameter $X$ of (1) is determined by the coefficients of the fifth degree equation. We first calculate the square root of the discriminant of the fifth degree equation in terms of the $z_1,\dotsc,z_5$. The transcendental part is to calculate the root $x$ of the icosahedral equation from the parameter $X$ by infinite processes. This is due to the hypergeometric series, as well as the transcendental solution of equation (2) by the binomial series. In the "Lectures on the Icosahedron," it has been shown, in particular, that all algebraic investigations which have been made for the purpose of solving the general equations of the fifth degree are reenactments of the aforementioned algebraic problem. The transcendental part of the task is only barely touched. It is clearly stated, however, what the connection is with the so-called 'solution of the equations of the fifth degree by elliptic functions.' I refer here to my other detailed explanations in the "Lectures on the Theory of Elliptic Modular Functions" \cite{KleinFrickeElliptic}, edited by Fricke and myself. There is a necessary connection between the fifth order transformation of the elliptic functions and the theory of the icosahedron. In (1), substituting $J$ (the absolute invariant of an elliptic modular function for $X$), the variable $x$ gets the meaning of the "principal modular function of the principal congruence group of the fifth degree."\footnote{Translator's Note: The terms "principal congruence [sub]group of the fifth degree" [of the modular group] and "principal modular function" are defined in \cite{KleinFrickeElliptic}.  See pages 388 and 591 in the original German books, or pages 323 and 475 in the English translation.} \\

All modes of relating the solution of the fifth degree equations to the elliptic functions are based on this fundamental theorem. In particular, $x$ can be represented by elliptic theta functions; it is a formula of principled simplicity, you have (if I may use Jacobi notation for the sake of brevity):
\begin{align}
    x = q^\frac{2}{5} \frac{\theta_1\lp \frac{2iK'\pi}{K},q^5 \rp}{\theta_1\lp \frac{iK'\pi}{K},q^5 \rp}
\end{align}

\noindent
\textit{However, the use of this formula to solve the icosahedral equation (or similar formulas for solving any resolvents of the Icosahedral equation) is just as much a detour as the solution of the pure equation (2) by logarithms:}
\begin{align}
    x = e^{ \frac{1}{n}\log(X) }
\end{align}

\noindent
You have to first calculate $\frac{K'}{K}$, respectively, by calculating $\log(X)$ from $X$ before applying formulas (3),(4). The meaning of the formulas for the solution is at most a practical one, namely if one has a logarithm table of elliptic periods $K,K'$. \textit{Thus, we can finally realize that the use of elliptic functions is not the essence of the theory of equations of the fifth degree.} This mode of expression via elliptic functions is only a residue of accidental historical development: the transformation theory of functions has given the first approach to establishing certain simple equations closely related to the icosahedral equation, namely the modular equations and multiplier equations for the fifth degree transformation. \\

So much for the introduction of the icosahedron into the theory of the fifth degree in general. I now have to limit myself to the algebraic side of the task. And here, above all else, I have to mention a fundamental proposition about the icosahedral substitutions, which becomes particularly important in what follows. One can pass from the icosahedral substitutions of the variable $x$ appearing in (1) to homogeneous substitution formulas (by replacing $x$ in the substitution formulas everywhere by $x_1:x_2$ and separating numerator and denominator in an appropriate manner). If one chooses the determinant of the resulting binary substitutions equal to 1, one has 120 binary substitutions; specifically, the identity substitution $x'=x$ corresponds to the two homogeneous substitution formulas
\begin{align}
    x_1' = x_1, x_2' = x_2 \quad \text{ and } \quad x_1' = -x_1, x_2' = -x_2.
\end{align}

\noindent
\textit{It is not possible in any way} (even if you change the value of the determinant), \textit{to assemble from such homogeneous substitutions a group which is isomorphic with the non-homogeneous substitution group (which contains fewer than 120 substitutions).} The surjective homomorphism from the substitution group of the $x_1:x_2$ to the substitution group of the $x$ therefore has a non-trivial kernel. This fundamental proposition, which is somewhat abstract, gives the algebraic theory of equations of the fifth degree its peculiar form, as we shall have to explain at once. Let us note in advance that it is not difficult to prove it. On pages 46 and 47 of my book on the icosahedron, it is traced back to the fact that the group of non-homogeneous icosahedral substitutions contains the Klein four-group and the corresponding proposition already applies to the Klein four-group. Let us take the following to be the simplest representation of the Klein four-group as (non-homogeneous) substitutions; it is given by
\begin{align}
    I: \xi' &= \xi   &II: \xi' &= - \xi   &III: \xi' &= \frac{1}{\xi}    &IV: \xi' = -\frac{1}{\xi}
\end{align}

\noindent
Here, $II$, $III$, and $IV$ are substitutions of order 2 and, at the same time,
\begin{align}
    II \cdot III \cdot IV = I.
\end{align}

\noindent
If one now wishes to create an isomorphic group of homogeneous substitutions, one certainly has a group with
\begin{align*}
    I': \xi_1' = \xi_1, \xi_2' = \xi_2
\end{align*}

\noindent
and replacing $II$, $III$, and $IV$ by
\begin{align*}
    II' &:  \xi_1' = \mp \xi_1, \xi_2' = \pm \xi_2 \\
    III' &: \xi_1' = \pm \xi_2, \xi_2' = \pm \xi_1 \\
    IV' &:  \xi_1' = \mp \xi_2, \xi_2' = \pm \xi_1 
\end{align*}

\noindent
(where, in the individual horizontal rows, the upper or lower signs are to be taken as desired). But, as one must also choose the signs here, \textit{the substitutions} $II'$, $III'$, \textit{ and } $IV'$ \textit{ each have determinant} $-1$ and thus it is impossible that
\begin{align*}
    II' \cdot III' \cdot IV' = I'
\end{align*}

\noindent
However, this contradicts (7) and no such isomorphism can exist. From now on, we will understand the homogeneous icosahedral substitutions as the 120 binary substitutions with determinant +1, corresponding to the 60 non-homogeneous substitutions of $x$.

\newpage
\noindent
I will now formulate the central problem for which we are responsible:
\begin{center}
    \textit{From the five independent variables $z_1,\dotsc,z_5$ (the roots of the equation of the fifth degree), one has to compose a function $x(z_1,\dotsc,z_5)$ which gives an isomorphism from the 60 icosahedral substitutions to the 60 even permutations of $z_1,\dotsc,z_5$.} 
\end{center}

\noindent
From our fundamental proposition, it immediately follows that there is no such \textit{rational} function of the five free variables (Lectures on the Theory of Elliptic Modular Functions \cite{KleinFrickeElliptic}, p.255). Namely, by dividing $x$ into coprime polynomials corresponding to the numerator and denominator; i.e. writing
\begin{align*}
    x(z_1,\dotsc,z_5) = \frac{\phi(z_1,\dotsc,z_5)}{\psi(z_1,\dotsc,z_5)}
\end{align*}

\noindent
with $\phi$ and $\psi$ coprime, the $\phi,\psi$ thus introduced would necessarily be homogeneously linear in the 60 permutations of the $z_1,\dotsc,z_5$. These homogeneous substitutions would correspond individually to the icosahedral substitutions of $x$. So, one would have a group of binary homogeneous substitutions that is isomorphic with the group of inhomogeneous icosahedral substitutions and such an isomorphism does not exist, as we have seen previously.\\

\textit{The required function $x(z_1,\dotsc,z_5)$ must therefore depend on its argument algebraically} \footnote{Translator's Note: In particular, this dependence is not rational.}. With this, we are led into the domain of those irrationalities of the theory of equations which I call \textit{accessory} in my Lectures (\cite{KleinLectures}, p. 158,159), because they are added as something new to the immediately existing irrationalities of the rational functions of the $z_1,\dotsc,z_5$. \footnote{Translator's Note: We refer to these in modern language as natural irrationalities.} The usual Galois theory of equations can only deal with the immediately existing irrationalities and not the accessory irrationalities, as usually happens when something new comes forward.\\

We do not know anything about the efficacy of these accessory irrationalities in general. Rather, we are dependent on tentative experiments in individual cases. Certainly, in the solution of any higher equation, the only allowable accessory irrationalities are those calculated from the symmetric functions of the roots (possibly the predetermined rational functions) by means of \textit{lower} equations. \footnote{Translator's Note: This is Klein's statement of resolvent degree.} In the equations of the fifth degree, which we treat here, the symmetric functions of the $z_1,\dotsc,z_5$ and its difference product (the square root of its discriminant) are known. 

\np
\textit{We can successfully construct (in many ways) one of the icosahedrally-dependent $x(z_1,\dotsc,z_5)$ as soon as we adjoin the square root of a suitable rational function of the $z_1,\dotsc,z_5$.} \footnote{The fifth root of unity $\epsilon = e^\frac{2\pi i}{5}$ occurs in the icosahedral substitutions and will be useful in the construction of a suitable function $x$. If we count the accessory irrationalities rigorously, then one has accessory irrationalities in the theory of equations from the beginning - namely, in the reduction of the cyclic equations to pure equations. \\} The two methods of solving equations of the fifth degree, which I give in my lectures, differ only by where they adjoin this accessory square root. In the first method, the accessory square root (transforming the fifth-degree equation by a Tschirnhaus transformation into a so-called fifth degree 'principal' equation - that is an equation in which the sum of the roots and the sum of the square roots vanishes) comes first. In the second method, we first take a step towards the icosahedral problem and then adjoin the accessory square root. As we said in the introduction, I give preference to the second method here, by introducing its individual steps in such a way that the whole approach can analogously be transferred to the sixth degree.\\

\quad Here, in numbered order, are the main considerations (of the second method):
\begin{enumerate}
    \item When $x_1:x_2$ undergo the 120 homogeneous binary icosahedral substitutions, the squares and the product
    \begin{align*}
        x_1^2, x_1x_2, x_2^2.
    \end{align*}
    
    \noindent
    undergo only 60 homogeneous ternary substitutions of determinant 1 (whose group is isomorphic to the 60 non-homogeneous icosahedral substitutions of $\frac{x_1}{x_2}$, and thus to the 60 even permutations of the five quantities $z_1,\dotsc,z_5$).
    \item The same is true, according to the general principles of invariant theory, of the coefficients of the quadratic binary form of $x_1:x_2$. In order to have an immediate connection to the style of Kronecker and Brioschi (also used in my Lectures), I shall hereby designate such a form as follows:
    \begin{align}
        A_1x_1^2 + 2A_0x_1x_2 - A_2 x_2^2
    \end{align}
    
    \noindent
    The $A_1,$ $2A_0,$ $A_2$ depend contragrediently on $x_1^2,$ $x_1x_2,$ $x_2^2$, according to the notions of invariant theory. \footnote{Translator's Note: Recall that given a group $G$ and linear $G$-representations $V$ and $W$, a $W$-valued contragredient of $V$ is a $G$-equivariant regular map $\AAA(V^\vee) \ra \AAA(W)$.}
    \item We readily conclude that it is possible to form (from any given five variables $z_1,\dotsc,z_5$) rational functions such that the alternating permutations of the $z_1,\dotsc,z_5$ also permute the $A_0,$ $A_1,$ $A_2$. In fact, in the work already mentioned in the introduction to Vol. 15 of the Math. Annalen, I have given a general approach which implies that whenever two given sets of variables (here the $z_1,\dotsc,z_5$ and the $A_0,A_1,A_2$) undergo isomorphic homogeneous linear substitutions, we can simply construct rational functions from which we can translate the first set of variables (like the $z_1,\dotsc,z_5$) into the second set (like the $A_0,A_1,A_2)$.
    \item We do not reproduce the general approach here (which would be unnecessarily lengthy), but instead give the abbreviated form relevant to our particular problem, which deals directly with the developments of Kronecker and Brioschi on equations of fifth degree \cite{Brioschi1858}. These are the following points:
    \begin{enumerate}
        \item We form six quadratic expressions in $x_1:x_2$:
        \begin{align}
            \sqrt{5}x_1x_2, \ \ \epsilon^\nu x_1^2 + x_1x_2 - \epsilon^{4\nu}x_2^2 \quad \text{ where } \epsilon = e^\frac{2\pi i}{5}, 0 \leq \nu \leq 4
        \end{align}
        
        \noindent
        which are permuted when $x_1:x_2$ undergo the icosahedral substitutions.
        \item Further, let $v(z_1,\dotsc,z_5)$ be a rational function of $z_1,\dotsc,z_5$, which remains invariant by the cyclic permutation of $z_1,\dotsc,z_5$ taken in natural order. We form the difference
        \begin{align*}
            v(z_1,\dotsc,z_5) - v(z_5,\dotsc,z_1)
        \end{align*}
        
        \noindent
        and square it. Then, we have a "metacyclic" function,\footnote{Translator's Note: This means that the corresponding stabilizer is a metacyclic group. Recall that a group is metacyclic if it is an extension of a cyclic group by a cyclic group. In this case, the corresponding stabilizer is $D_{10}$, the dihedral group with 10 elements, which is indeed metacyclic. \\} which should be called $u_\infty^2$, while the five other values that result from it by the alternating permutations of $z_1,\dotsc,z_5$ may be labeled $u_\nu^2$ in a proper order $(\nu = 0,1,2,3,4)$. One can then choose the signs of the $u_\infty,u_0,\dotsc,u_4$ \textit{such that the alternating permutations of the $z_1,\dotsc,z_5$ permute the}
        \begin{align}
            u_\infty, u_0, \dotsc, u_4
        \end{align}
        
        \noindent
        \textit{isomorphically by the corresponding icosahedral substitutions of $x_1,x_2$; namely, they undergo the same sign changes as (9)}.
        \item We conclude that the following form
        \begin{align}
            \Omega(z_1,\dotsc,z_5|x_1:x_2) = \sqrt{5}u_\infty x_1x_2 + \SL_{\nu=0}^4 u_\nu \lp \epsilon^\nu x_1^2 + x_1x_2 - \epsilon^{4\nu}x_2^2 \rp
        \end{align}
        
        \noindent
        remains invariant if one simultaneously applies an alternating permutation to the $z_1,\dotsc,z_5$ and the corresponding icosahedral substitution to the $x_1:x_2$.
        \item We now put, in accordance with (8):
        \begin{align*}
            \Omega(z_1,\dotsc,z_5|x_1,x_2) = A_1x_1^2 + 2A_0x_1x_2 - A_2x_2^2
        \end{align*}
        
        \noindent
        and find by comparison
        \begin{align}
            \begin{cases}
                2A_0 &= \sqrt{5}u_\infty + \SL_{\nu=0}^4 u_\nu\\
                A_1 &= \SL_{\nu=0}^4 \epsilon^\nu u_\nu\\
                A_2 &= \SL_{\nu=0}^4 \epsilon^{4\nu} u_\nu
            \end{cases}
        \end{align}
    \end{enumerate}
    \textit{Thus, we have constructed from $z_1,\dotsc,z_5$ the quantities $A_0, A_1, A_2$ which are permuted in the desired way when the $z_1,\dotsc,z_5$ undergo an alternating permutation.}
    \item We refer to the above result by saying that we have assigned a covariant quadratic binary form (8) to the $z_1,\dotsc,z_5$. \footnote{Translator's Note: Given a group $G$ and linear $G$-representations $V$ and $W$, a $W$-valued covariant of $V$ is a $G$-equivariant regular map $\AAA(V) \ra \AAA(W)$. \\} The discriminant of (8) is
    \begin{align}
        A = A_0^2 + A_1A_2,
    \end{align}
    
    \noindent
    which is a binary function of $z_1,\dotsc,z_5$ that is invariant under alternating permutations of $z_1,\dotsc,z_5$; it is thus a rational function of the coefficients of the given fifth degree equation and the square root of its discriminant. The goal is not to assign to $z_1,\dotsc,z_5$ a covariant binary quadratic form or a "pair of points" of the binary form
    \begin{align}
        A_1x_1^2 + 2A_0x_1x_2 - A_2x_2^2 = 0,
    \end{align}
    
    \noindent
    but to assign a quotient $\frac{x_1}{x_2}$, i.e. a point. We do this in the simplest way by solving the quadratic equation (14)  and accordingly writing
    \begin{align}
        \frac{x_1}{x_2} = x = \frac{-A_0 + \sqrt{A_0^2+A_1A_2}}{A_1}.
    \end{align}
    \item \textit{Thus, we have solved our central task: to compose such an $x$ from the $z_1,\dotsc,z_5$ which undergoes the icosahedral substitutions corresponding to the alternating permutations of the $z_1,\dotsc,z_5$}. Note that the $A_0,A_1,A_2$ of item 4d are rational functions of $z_1,\dotsc,z_5$ and their construction uses only one irrationality - the fifth root of unity $\epsilon$. According to item 5, we see that the expression of $A_0,A_1,A_2$ under the square root\footnote{Translator's note: e.g. the $A_0^2+A_1A_2$ in equation (15)} is invariant under alternating permutations of the $z_1,\dotsc,z_5$. \textit{We have thus achieved the goal with the aid of such accessory irrationalities, which in the theory of the fifth degree, will suitably be called lower irrationalities.}
    \item We now further investigate the parameter $X$ of the icosahedral equation, which satisfies our $x$ (15) as a function of the coefficients of the equation of the fifth degree whose roots are the $z_1,\dotsc,z_5$, respectively. We will say that the equation is solved if we calculate the square root of the discriminant and how the $z_1,\dotsc,z_5$ are rationally represented by $x$ from the coefficients [of the equation of the fifth degree] and the adjoined square root [i.e. we do not concern ourselves with the transcendental portion].
    \item In summary, let us emphasize why one can justifiably speak of such a solution of equations of the fifth degree. Not only is there a sequence of steps that could be numerically traversed (in the given case) so that one actually obtains the numerical values of $z_1,\dotsc,z_5$, but it is also a full theoretical insight into the internal nature of the problem of solution. \footnote{Translator's note: By "problem of solution", Klein means the problem of solving generic polynomials. \\} After all, the $z_1,\dotsc,z_5$ are the different branches of a finite-valued algebraic function, which depends on the coefficients of the fifth degree equation and is initially of a very confusing design. \textit{These five branches $z_1,\dotsc,z_5$ are defined over the field of rationality determined by the icosahedral irrationality. \footnote{Translator's Note: For more on fields of rationality, see Ackerman's English translation of "Development of Mathematics in the 19th Century" by Felix Klein. In particular, see Chapter VII - Deeper Insight into the Nature of Algebraic Varieties and Structures, Section 3: The Theory of Algebraic Integers and Its Interaction with the Theory of Algebraic Functions, p.312-314.} The field of rationality is given by the coefficients of the fifth degree equation, the square root of the discriminant, and the adjoined accessory irrationalities. The icosahedral irrationality is of a more transparent construction and is a higher irrationality that depends only on a single parameter from the field of rationality.}
    %%%
    %%% This was my previous translation, before talking with Yao
    %%%
    %%% \textit{These five branches $z_1,\dotsc,z_5$ are defined over one of the extension fields of the base field of the fifth degree equation obtained by the square root of their discriminant and the multiple accessory irrationalities. Note that the multiple accessory irrationalities can be replaced by a single, more transparent irrationality that depends only on a parameter determined by the icosahedral irrationality.} \todo{have Jesse verify}
\end{enumerate}

I would like to take up at this point a more personal remark about the relationship between my papers on the equations of the fifth degree and those of Kronecker - as you, dear colleague, are in the position of being able to see the manuscripts of Kronecker and thus can complete my information in an authentic way. As is well known, Kronecker and Brioschi \cite{Brioschi1858} used the same quantities $A_0,A_1,A_2$ in their first papers on equations of the fifth degree (from 1858) which I quoted earlier (in list item 4b); they then constructed the sixth-degree equation satisfying $\zeta = 5A_0^2$, and which Brioschi calls a "Jacobi equation" because of its close connection with certain equations established by Jacobi for the transformation of elliptic functions; finally, they state that by adjoining a root, one can arrive at an equation with only one parameter. This square root defines an accessory irrationality equivalent to that used in formula (15). Furthermore, Kronecker \cite{Kronecker1861} set up the fundamental theorem, which I designate as Kronecker's theorem in my Lectures. Moreover, the exposition and proof of Kronecker's theorem are the crowning achievement of my Lectures; the theorem is that \textit{it is impossible to form a resolvent to the general equation of the fifth degree with only one parameter and without resorting to accessory irrationalities}. As in the 12th volume of the Math. Annalen \cite{Klein1877}, I prove this proposition by invoking the property of the icosahedral group discussed above, namely the doubling (at least) of its substitutions in the transition to homogeneous substitutions. My first proof, which I gave in 1877 in the reports of the Erlangen physical-medical school (meeting of January 13), was considerably more complicated. Twenty-four years ago (Easter, 1881), I had the opportunity to talk in detail with Kronecker about these things. It turned out that in his investigations, Kronecker was unaware of the icosahedral substitutions to which he had come so close and, accordingly, did not have sufficient proof for his main claim! I think this is a very remarkable fact, but also very common, for it confirms in a particularly interesting case what Gauss so often emphasizes: that the discovery of the most important mathematical theorems is more a matter of intuition than of deduction, and the production of the proof is a very different business from the discovery of the theorems. I did not return to the subject later with Kronecker, but some years ago, I heard that after the publication of my Lectures in a college, Kronecker has commented on the solution of the fifth degree equations and the theory of the icosahedron. I would be very interested (as certainly would other mathematicians) in finding out what may be contained in Kronecker's papers on these matters, and I would like to ask you to review the relevant material and publish it soon. \\

A new proof of Kronecker's theorem has been given by Mr. Gordan in Volume 29 of the Math. Annalen \cite{Gordan1887}. \footnote{Compare this with the presentation of Gordan's proof in the textbook of Weber and Netto's "Algebra" \\} It is easier to read than mine, in that it does not refer to an explicit knowledge of the icosahedral substitutions anywhere. Nevertheless, as I shall point out, it is most closely related to the basic idea of my proof. Following a development by Mr. L\"uroth \cite{Luroth}, we both use a proposition which can be formulated as follows: \footnote{Translator's Note: The theorem of L\"uroth that Klein and Kronecker use can be stated as follows - Let $\kk$ be an algebraically closed field of characteristic 0. Then, any unirational function field of transcendence degree 1 is isomorphic to $\kk(t)$, where $t$ is transcendental over $\kk$. }\\

Suppose an equation of $n^{th}$ degree, whose roots are the independent variables $z_1,\dotsc,z_n$, has a rational resolvent with only one parameter. Then, the $n^{th}$ degree equation must have a rational function $x$ of the $z_1,\dotsc,z_n$ such that when the $z_1,\dotsc,z_n$ are permuted by the Galois group, $x$ undergoes a linear transformation. We have the same understanding that, when changing to homogeneous coordinates $x_1:x_2$, we must go from our group of linear substitutions to an isomorphic groups of homogeneous linear substitutions. On the other hand, of course, this group must be isomorphic to the Galois group of the equation modulo a specified subgroup. \\

Now, on p.44-47 of my Lectures, I have given the proposition that only the following groups of linear substitutions of a single variable can be converted isomorphically to their corresponding binary forms:
\begin{enumerate}
    \item The cyclic groups
    \item The dihedral groups of odd $n$. 
\end{enumerate}

\noindent
It follows that \textit{an equation of the $n^{th}$ degree (whose roots are the independent variables $z_1,\dotsc,z_n$) only admits a rational resolvent with only one parameter (which can be immediately transformed into a pure equation, or dihedral equation of odd $n$) if its Galois group is (isomorphic to) a normal subgroup of a cyclic group or a dihedral group of odd $n$.} \footnote{Translator's Note: In modern language, this is the classification of groups of essential dimension 1.} An associated resolvent with only one parameter can then be set up immediately according to the principles in the Math. Annalen Vol 15 \cite{Klein1879}. The general statement above is included in both Gordan's proof and my proof of Kronecker's theorem. Indeed, my proof is done by pointing out that the group of a fifth degree equation with the square root of the discriminant adjoined is simple; however, it is isomorphic to the group of linear substitutions of the icosahedron and thus the previous theorem yields the claim. Gordan's proof, on the other hand (if I understand it correctly), uses the obvious fact that the group in question, like every group, contains itself as a normal subgroup. The quotient of the group by itself is isomorphic to the identity group. And the identical substitution falls under the premise of our theorem. Thus, there are in fact resolvents with one parameter, but they are completely useless for the solution of the equations of the fifth degree! Namely, there are linear resolvents whose square root is a function of the $z_1,\dotsc,z_5$ that is invariant under the alternating permutations of the $z_1,\dotsc,z_5$. But there are no other (rational) resolvents with just one parameter, or better: every rational resolvent of our fifth degree equation with just one parameter is linear and therefore useless. \\

So much for the icosahedral substitutions and the solution of fifth degree equations mediated by them. Instead of the "unitary" substitutions $x' = e^\frac{2\pi ik}{n}$, which link the roots of a \textit{pure} equation to one another, "binary" linear substitutions of two homogeneous variables $x_1:x_2$ have entered. At the same time, the way to new generalizations has opened up - one simply has to use groups of linear substitutions of several homogeneous variables! I cannot possibly repeat here the reflections which I gave in this regard first in the 15th volume of the Math. Annalen \cite{Klein1879} or recall the elaborations which have later been concluded. It suffices to refer to Weber's textbook \cite{Weber} and to the already cited encyclopedia article by Wiman \cite{Wiman1898}. \footnote{A first clear introduction is also the Lecture IX of my Evanston Colloquium, held at the World's Fair in Chicago (Macmillan, New York, 1894)} In this sequence, we consider \textit{an equation of sixth degree along with the square root of its discriminant}, whose Galois group consists of the 360 alternating permutations of the roots $z_1,\dotsc,z_6$. It will be necessary to use the smallest number of homogeneous $x_1:\dotsc:x_{\mu}$ for which there exists a surjective homomorphism from the group of linear substitutions of the $x_1:\cdots:x_\mu$ to the group of the 360 alternating permutations. If this surjective homomorphism were to prove to be an isomorphism, we would be able to write down rational functions of $z_1,\dotsc,z_6$ (according to the prescriptions of \cite{Klein1879}) which are linear in the 360 alternating permutations of $z_1,\dotsc,z_6$ and yield $x_1:\dotsc:x_\mu$ as a result. However, it turns out that here, as in the equations of the fifth degree, the homomorphism must have a non-trivial kernel, so that we are asked the question whether (or respectively, how) we can get by with the help of lower accessory irrationalities. \\

My first approach to the formulation of this question is in Volume 28 of the Math. Annalen (1886, On the Theory of General Equations of the Sixth and Seventh Degrees)\footnote{Translator's Note: Klein gives the correct volume for the article, but the incorrect year. This article is from December 1887 in Issue 4 of Volume 28. \\} \cite{Klein1887}. At the time, it seemed like the preliminary work of Mr. C. Jordan would not allow a ternary group of linear substitutions with the requisite surjective homomorphism to the group of the 360 alternating permutations on six letters; such a group was first discovered by Mr. Valentiner in 1889 (Volume 6 of Series V of the Danish Academy's papers: The Definitions of the Final Transformation Groups) \cite{Valentiner} and examined by structural and related fundamental invariants for the first time by Mr. Wiman in 1895 \footnote{Translator's Note: While Klein gives the correct volume, this work was actually published in December of 1896. \\} (Math. Annalen 47: About the Simple Group of 360 Plane Collineations) \cite{Wiman1896}. At that time, I constructed a group of quarternary collineations that is isomorphic to the group for the general equation of degree six - and also for the general equation of degree seven - and showed that the original problems rest on the corresponding problems for the groups of quarternary collineations, which can be obtained from the general equations of degree six (respectively, seven) with at most two accessory square roots. \footnote{The group which I put forward for the sixth degree equation contains as many as 720 collineations, so that it is not necessary to use all of them; adjoin the square root of the discriminant of the sixth degree equation beforehand. In contrast, the group corresponding to the equations of the seventh degree contains only $\frac{7!}{2} = 2520$ collineations. \\} \\

\textit{As far as the equations of the sixth degree are concerned (to which we confine ourselves here), this approach is currently unnecessary, because of the discovery of the Valentiner group.} \footnote{For the equations of seventh degree, the quaternary approach persists; but it is impossible to pursue the interesting questions in this text. \\} I note this expressly, because this is where Mr. Lachtin makes an unnecessary detour (as mentioned at the beginning of this letter). In order to connect the equation of degree six with the Valentiner group, Mr. Lachtin goes through the development given in Volume 28 of the Math. Annalen. \footnote{Translator's Note: Lachtin only has two articles in the Mathematische Annalen. The one Klein is referring to from September 1898 in Issue 3 of Volume 51. Lachtin's other article is on the septic and is from September 1902 in Issue 3 of Volume 56. \\} \cite{Lachtin1898} \\

This is not uninteresting, \footnote{Mr. Lachtin notes that in the quaternary group, the degree two surfaces in space interchange in much the same way as the degree three curves of the plane. From here, as is noted in passing, it is possible without any great difficulty to arrive at the same $\Sigma$, which I communicate below under (19). One only has to keep in mind that the roots $z_1,\dotsc,z_6$ of the sixth degree equation, and also their squares $z_1^2,\dotsc,z_6^2$ define a linear complex in space according to the developments of Volume 28, and that these two complexes together with the "unitary complex" introduced there determine a degree two surface through their common lines. No accessory irrationality occurs here. It is then in no way necessary, in the transition from space to plane, to refer to the comparatively complicated formulas, as Mr. Lachtin does, by which in Volume 28 I have assigned a point of space to the roots $z_1,\dotsc,z_6$, so also in this regard, the approach of Mr. Lachtin can be shortened.} but is by no means necessary for what we do next. The transition from the equations of the sixth degree to the Valentiner group (as I suggested in my Roman note of 1899) and which I now wish to introduce in more detail, does not require any reference to the quarternary substitution group. For the sake of clarity, I shall again divide the considerations in question into a numbered list below which makes clear the analogy with the train of thought followed for the equations of the fifth degree:

\begin{enumerate}
    \item The task is to compose three functions $x_1,x_2,x_3$ from the six free variables $z_1,\dotsc,z_6$ such that the homogeneous $x_1:x_2:x_3$ undergo the corresponding collineations of the Valentiner group when the $z_1,\dotsc,z_6$ undergo one of the 360 alternating permutations.
    \item Now, Mr. Wiman has already noticed that the number of collineations from the Valentiner group at least triples when one goes from substitutions of the plane to the corresponding ternary linear substitutions. Therefore, it is not possible for the required $x_1,x_2,x_3$ to be rational functions of $z_1,\dotsc,z_6$.
    \item From now on, we want to fix the homogeneous linear Valentiner substitutions so that their determinant is always 1. Thus, there are exactly $3(360) = 1080$ of them and the three substitutions corresponding to the identity substitution [of the plane] are:
    \begin{align}
        \begin{cases}
            x_1' = j^\nu x_1, \\
            x_2' = j^\nu x_2, \\
            x_3' = j^\nu x_3, 
        \end{cases}
        \quad j = e^\frac{2\pi i}{3}, \ \  0 \leq \nu \leq 2
    \end{align}
    \item We now note that in these 1080 homogeneous substitutions, the ten degree three terms coming from the $x_1,x_2,x_3$ are:
    \begin{align*}
        x_1^3,  \ x_1^2x_2, \dotsc
    \end{align*}
    
    \noindent
    which only undergo 360 homogeneous linear substitutions (and whose group will be isomorphic with the group of alternating permutations of $z_1,\dotsc,z_6$).
    \item We now consider any cubic ternary form
    \begin{align*}
        a_{1,1,1}x_1^3 + 3a_{1,1,2} x_1^2x_2 + \cdots
    \end{align*}
    
    \noindent
    (which, if set equal to 0, represents a "degree three curve" in the plane of $x_1,x_2,x_3$). The coefficients $a_{1,1,1}, 3a_{1,1,2}, \dotsc$ for any homogeneous linear substitutions of $x_1:x_2:x_3$ are related to the $x_1^3, x_1^2x_2, \dotsc$ contravariantly. Thus, in the substitutions of the Valentiner group, they also undergo exactly 360 homogeneous linear substitutions, which can be uniquely identified with the 360 alternating permutations of $z_1,\dots,z_6$. 
    \item We readily conclude that it is possible to form ten rational functions of the free variables $z_1,\dotsc,z_6$ denoted by:
    \begin{align*}
        \phi_{1,1,1}, \phi_{1,1,2}, \dotsc 
    \end{align*}
    
    \noindent
    which, in the case of the alternating substitutions of the $z_1,\dotsc,z_6$, are substituted just as the
    \begin{align*}
        a_{1,1,1}, a_{1,1,2}, \dotsc
    \end{align*}
    
    \noindent
    are by the corresponding substitutions of the Valentiner group; i.e., the roots $z_1,\dotsc,z_6$ rationally and covariantly assign a degree three curve.
    \item To put it another way, \textit{one can construct} (in many different ways and without the use of accessory irrationalities \footnote{Apart, of course, from the numerical irrationalities that occur in the substitutions of the Valentiner group. These are (in accordance with the following formulas (18), etc.) the square roots $\sqrt{-3}$ and $\sqrt{5}$.}), \textit{a cubic form depending on the $z_1,\dotsc,z_6$ and $x_1,x_2,x_3$}
    \begin{align}
        \Omega(z_1,\dotsc,z_6|x_1,x_2,x_3) = \phi_{1,1,1}x_1^3 + 3\phi_{1,1,2}x_1^2x_2 + \cdots 
    \end{align}
    
    \noindent
    \textit{which remains invariant if one simultaneously performs an alternating permutation on the $z_1,\dotsc,z_6$ and its corresponding Valentiner substitution on the $x_1,x_2,x_3$.}
    \item As far as the actual construction of such a form $\Omega$ is concerned, I do not give the general and extensive process which I provided in \cite{Klein1879}, but develop it in an abbreviated manner that emerged from my correspondence with Mr. Gordan (last winter), just as with equations of the fifth degree. One has to combine the following relationships:
    \begin{enumerate}
        \item The 360 collineations of the Valentiner group play an important role in two systems of six conic sections, as Mr. Wiman proved first. \textit{The six conic sections of each of the two systems are permuted by the corresponding 360 collineations in 360 ways}.
        \item The equations of these $2 \cdot 6 = 12$ conic sections were first proposed by Mr. Gerbaldi \cite{Gerbaldi} (Rendiconti del Circolo Matematico di Palermo, t. XII, 1898: "Sul gruppo semplice di 360 collineazioni piane, I; already published in 1882 in the Atti di Torino, Vol. XV, p. 358 ff., Note: "Sui gruppi di sei coniche in involuzione")\footnote{Translator's Note: The modern reference for this article is the one appearing in the Mathematische Annalen, which is what we give. Thus, we have left the two previous journal references Klein gave in the main text.}. Here we equate the corresponding ternary quadratic forms of determinant 1 by restricting to the one system of six conic sections, according to the procedure of Mr. Gordan.  We can then write: \footnote{Translator's note: As in equation (16), $j$ is the primitive 3rd root of unit $e^{\frac{2\pi i}{3}}$.}
        \begin{align}
            \begin{cases}
                k_1 &= x_1^2 + jx_2^2 + j^2x_3^2\\
                k_2 &= x_1^2 + j^2x_2^2 + jx_3^2\\
                k_3 &= \alpha\lp x_1^2 + x_2^2 + x_3^2 \rp + \beta\lp x_2x_3 + x_3x_1 + x_1x_2 \rp \quad \\
                k_4 &= \alpha\lp x_1^2 + x_2^2 + x_3^2 \rp + \beta\lp x_2x_3 - x_3x_1 - x_1x_2 \rp  \quad \\
                k_5 &= \alpha\lp x_1^2 + x_2^2 + x_3^2 \rp + \beta\lp -x_2x_3 + x_3x_1 - x_1x_2 \rp  \quad \\
                k_6 &= \alpha\lp x_1^2 + x_2^2 + x_3^2 \rp + \beta\lp -x_2x_3 - x_3x_1 + x_1x_2 \rp\\
            \end{cases}
        \end{align}
        
        \noindent
        where $\alpha = \frac{1-\sqrt{-15}}{8}$ and $\beta = \frac{ -3 + \sqrt{-15} }{4}$.
        \item The $k_1,\dotsc,k_6$ are determined, up to a third root of unity, by the requirement that their determinant be 1. In fact, the 1080 Valentiner substitutions permute the $k_1,\dotsc,k_6$ by multiplication by certain third roots of unity.
        \item We want, now, from any three of the $k$:
        \begin{align*}
            k', \ k'', \ k'''
        \end{align*}
        
        \noindent
        one whose coefficients form a trilinear covariant and a similar invariant. For the former, we choose the functional determinant $|k'k''k'''|$, which changes its sign when two of the $k',k'',k'''$ are exchanged. As an invariant, we take a symmetric combination of the coefficients of the $k',k'',k'''$ (namely the expression used in the development of the coefficient-determinant of the form $\lambda'k' + \lambda''k'' + \lambda'''k'''$ in which $\lambda'\lambda''\lambda'''$ appears). I will temporarily call it $(k'k''k''')$ here; this is a simple numerical quantity in the present case.
        \item For all possible triples $k',k'',k'''$, we now form the quotient
        \begin{align*}
            \frac{|k'k''k'''|}{(k'k''k''')}.
        \end{align*}
        
        \noindent
        One can show \textit{that the $\binom{6}{3} = 20$ quotients obtained above undergo the same sign changes from the 1080 substitutions of the Valentiner groups that the 20 difference products}
        \begin{align*}
            (z''-z''')(z'''-z')(z'-z'')
        \end{align*}
        
        \noindent
        \textit{undergo from the corresponding alternating permutations of the $z_1,\dotsc,z_6$.}
        \item Therefore, the sum of all triples
        \begin{align}
            \sum (z''-z''')(z'''-z')(z'-z'') \cdot \frac{|k'k''k'''|}{(k'k''k''')}
        \end{align}
        
        \noindent
        is a simple example of a form
        \begin{align*}
            \Omega\lp z_1,\dotsc,z_6,x_1,x_2,x_3 \rp
        \end{align*}
        
        \noindent
        as we were looking for in item number 7 above.
        \item More general examples (which we do not need in the following) are obtained by substituting the determinant
        \begin{align*}
            \lav
            \begin{matrix}
                z'^\alpha &z''^\alpha &z'''^\alpha\\
                z'^\beta  &z''^\beta  &z'''^\beta\\
                z'^\gamma &z''^\gamma &z'''^\gamma
            \end{matrix}
            \rav
        \end{align*}
        for the difference product of $z',z'',z'''$ in (19).
        \item Let us now assign the sum (19) to the successive terms $x_1^3,x_1^2x_2,\dotsc$ by writing, as in formula (17): 
        \begin{align}
            \sum (z''-z''')(z'''-z')(z'-z'') \cdot \frac{|k'k''k'''|}{(k'k''k''')} = \phi_{1,1,1}x_1^3 + 3\phi_{1,1,2}x_1^2x_2 + \cdots , 
        \end{align}
        
        \noindent
        so that the $\phi_{1,1,1}, \phi_{1,1,2}, \dotsc,$ are just such the rational functions of $z_1,\dotsc,z_6$ that we looked for in item number 6 above.
    \end{enumerate}
    \item The degree three curve 
    \begin{align}
        \sum (z''-z''')(z'''-z')(z'-z'') \cdot \frac{|k'k''k'''|}{(k'k''k''')} = 0 
    \end{align}
    
    \noindent
    (whose coefficients depend rationally on the $z_1,\dotsc,z_6$) covariantly assigns a point $x_1:x_2:x_3$ using accessory irrationalities that are as low as possible. 
    \item The theory of degree three plane curves offers various possibilities. For the sake of brevity, as I did in my Roman note, I want to choose an inflection point in third-order curve here.
    \item According to the well-known theory of Hesse, the determination of such an inflection point requires only square roots and cube roots; for the solution of the equations of the sixth degree, these are indeed lower irrationalities. The details [of the process of determining an inflection point] are not discussed here.
    \item On the other hand, the inflection point is certainly connected in a covariant manner with the degree three curve: if any collineation is performed on a curve with a chosen inflection point $P$, then $P$ is taken to an inflection point $P'$  on the new curve and each of the nine inflection points can be realized in this manner. In particular, this applies to the 360 collineations of the Valentiner group.  
    \item We now think of the coordinates $x_1:x_2:x_3$ of our chosen inflection point, instead of the coefficients $\phi_{1,1,1},\phi_{1,1,2},\dotsc$ of the curve of third order (whose values come from (20) and the $z_1,\dotsc,z_6$).
    \item If the $z_1,\dotsc,z_6$ undergo an alternating permutation, the $x_1:x_2:x_3$ undergo the corresponding collineation of the Valentiner group.\\
    \\
    We conclude that the rational functions of $z_1,\dotsc,z_6$ that remain invariant after the reductions in the expressions of $x_1:x_2:x_3$ from the occurring square roots and cubic roots also remain invariant under the alternating permutations of the $z_1,\dotsc,z_6$. Thus, they can be represented as rational functions of the coefficients of the presented sixth degree equation and the square root of their discriminant. 
    \item Therefore, we will are justified in designating the irrationalities required in the calculation of the inflection point as \textit{lower} accessory irrationalities. 
    \item By computing the coordinates $x_1:x_2:x_3$ of an inflection point of our $C_3$ [our degree three curve], we have accomplished the goal: \textit{we form the functions $x_1, x_2, x_3$ from the free variables $z_1,\dotsc,z_6$ and lower accessory irrationalities such that when the the $z_1,\dotsc,z_6$ undergo an alternating permutation, the $x_1, x_2, x_3$ undergo the corresponding collineation of the Valentiner group.}
\end{enumerate}

This is the explanation of the particular content of my Roman note, which I thought to give here. \footnote{The final sentence of the note has become incomprehensible when printed in the Rendiconti of the Accademia dei Lincei by a strange change. It should read "And with the aid of accessory irrationalities, which we usually regard as elementary, we come to the goal." Instead, what is printed is "And so, with the aid of ancillary irrationalities, we come to the goal, which usually regarded as elementary."} \\

One might wish for a closer examination of the inference used in list item number 14. The simplest thing would be to calculate all the known equations leading to the determination of an inflection point of the degree three curve (21) and thus actually confirm the assertion. Moreover, as Mr. Gordan remarks, the whole of the conclusion can be dealt  with in the following way. Just consider the ninth degree equation, which satisfies the nine values that an absolute invariant of the Valentiner group (e.g. the $\nu$ to be named immediately) assumes in the nine inflection points! This equation must be the same for all 360 third-order curves (which results from the substitutions of the Valentiner group and thus by the alternating  substitutions of the two). After disposing of indifferent factors,  the coefficients are rational functions of the $z_1,\dotsc,z_6$ that are invariant under the alternating permutations of the $z_1,\dotsc,z_6$. The affect of this ninth degree equation can be none other than that of the original inflection point equation. It is thus solved by square roots and cubic roots of rational functions of $z_1,\dotsc,z_6$ that are invariant under the alternating permutations of $z_1,\dotsc,z_6$.\\

Now, if we adjoin one of the resulting nine values of our absolute invariant ($\nu$), then it and the equation (21) of the third order curve, (respectively, of the equation of its Hessian curve), becomes the corresponding single inflection point $x_1:x_2:x_3$ and is calculated rationally. Consequently, the assertion of list item number 14 concerning the irrationalities required in the calculation of the inflection point, is self-evident. \\

The further treatment of equations of the sixth degree will have to be done, in any case, by calculating the absolute invariants of the Valentiner group of the selected inflection point of our third order curve. According to Mr. Wiman, the Valentiner group has three lowest invariants:
\begin{align}
    F, \ H, \ \Phi
\end{align}

\noindent
of degrees 6, 12, and 30 in the $x_1,x_2,x_3$. From them, the two fundamental absolute invariants come together, which I call $v$ and $w$ here in connection with the work of Mr. Lachtin to be mentioned immediately:
\begin{align}
    v = \frac{\Phi}{F^3}, \quad w = \frac{H}{F^2}.
\end{align}

\noindent
If we enter the coordinates of our point of inflection for $x_1:x_2:x_3$, then the $v, w$ are rational functions of the coefficients of the sixth degree equation and the square root of the discriminant (respectively, the occasionally introduced accessory irrationalities). \textit{The} Normalproblem \footnote{In German, this was indeed one word - \textit{Normalproblem}. \\} \textit{of solving equations of the sixth degree is thus the reduction to calculating $x_1:x_2:x_3$ from the known $v,w$}. As just stated, this is now a problem with two arbitrary parameters and is distinguished by the fact that all of its 360 solutions $x_1:x_2:x_3$ can be determined from a given solution by the 360 collineations of the Valentiner group. We do not currently have a method to reduce the number of parameters to one by means of further lower irrationalities. For example, if we try to assign a point $x_1':x_2':x_3'$ to the degree six curve $F = 0$ in a covariant manner and thus set the following (instead of the Normalproblem (23) to the inflection point $x_1:x_2:x_3$)
\begin{align}
    F' = 0, \quad t' = \frac{\Phi'^2}{H'^5},
\end{align}

\noindent
in the usual approach (intersection of the curve $F=0$ with a straight line covariantly dependent on the point $x_1:x_2:x_3:$), one encounters an auxiliary equation which itself is of the sixth degree! \\

For the sake of completeness, we finally ask for a reverse form, i.e. to rationally calculate the quantities assumed to be known [$v,w$] from the roots $z_1,\dotsc,z_6$ of the presented sixth order equation and a single solution system $x_1:x_2:x_3$ of (23). \textit{Thus, we have completely sketched the algebraic part of the solution of the equations of the sixth degree.} \\

\textit{The transcendental part will require infinite processes to actually compute the $x_1:x_2:x_3$ from the equations (23).} A first approach to this is done by Mr. Lachtin in a voluminuous work, which was first published in Russian (1901) in the 22nd volume of the Moscow Mathematics Collection and then in 1902 in the German edition of Volume 56 of the Math. Annalen \cite{Lachtin1902}. \footnote{"The differential resolvent of an algebraic equation of the sixth degree of a general kind." (Math. Ann., Vol. 56, p. 445-481.) \\} Writing
\begin{align}
    y_1 = \frac{x_1}{ \sqrt[6]{F} }, \quad y_2 = \frac{x_2}{ \sqrt[6]{F} }, \quad y_3 = \frac{x_3}{ \sqrt[6]{F} },
\end{align}

\noindent
the $y_1,y_2,y_3$ are a system of solutions to the three simultaneous linear partial differential equations expressing the second derivatives
\begin{align*}
    \frac{\p^2 y}{\p v^2}, \quad \frac{\p^2 y}{\p v \p w}, \quad \frac{\p^2 y}{\p w^2}
\end{align*}

\noindent
linearly in the $\lp \frac{\p y}{\p v} \text{ and } \frac{\p y}{\p w} \rp$ and the $y$.  Mr. Lachtin has shown that the coefficients of these equations are rational functions of the absolute invariants $v,w$, which do not exceed certain definable degrees. However, he did not calculate the numerical coefficients of these polynomials. The remaining gap is now being filled by the work of Mr. Gordan, which I referred to in the beginning of this letter. \footnote{We believe that Klein is referring to the work that became \cite{Gordan}} \textit{In fact, Mr. Gordan succeeded in making explicit the partial differential equations in question}. It is thus possible to develop the $y_1, y_2, y_3$ in powers of $v$ and $w$, or even in any series of linear functions of $v$ and linear functions of $w$; thus, it is no longer an issue to determine the regions in which the various series thus formed converge - \textit{in other words, we can solve the transcendental problem in a direct way.} \\

Again, for the sake of completeness, it must be added that the special problem presented by equation (24) is already discussed in detail in terms of function theory. In 1896, at the Frankfurter Scientific Congress, Mr. Fricke \footnote{Translator's Note: For more information, see \cite{Fricke} } dealt with the decomposition of the Riemann surface (of genus 10) corresponding to the Valentiner group into fundamental domains and a closed relation of the same with the decomposition of the half plane in the semicircular triangles from the angles
\begin{align*}
    \frac{\pi}{2}, \quad \frac{\pi}{4}, \quad \frac{\pi}{5}
\end{align*}

\noindent
Mr. Lachtin then confirms this information in \cite{Lachtin1898} and established the third-order linear differential equation, for which - in the case of equations (24) - the parameter $t$ is satisfied by the variables $x_\nu$ multiplied by a suitable factor.\\

I am at the end of my presentation. I hope the analogy of the proposed sixth degree equations with the solution of the equations of the fifth degree by the icosahedral equation appears convincing. A finer examination of the details given by Mr. Gordan and myself for the equations of the fifth degree, as well as a geometric presentation, are given in my "Lectures on the Icosahedron." \\

G\"ottingen, March 22, 1905.

\np
\section{Ending Footnote}
Mr. Gordan has been able to devote only one introductory essay to the questions raised above \cite{Gordan1909}. There he makes a substantial simplification of the necessary accessory irrationality of $x_1,x_2,x_3$. Instead of the degree three curve of the $x_1,x_2,x_3$-plane, which I rationally assigned to the value system $z_1,\dotsc,z_6$, it uses a (1,1)-connection; i.e. a bilinear form in the $x$ and $u$ (whose coefficients must be assumed to be whole rational functions of $z_1,\dotsc,z_6$ and that the form remains invariant under the 360 permutations of the $z_1,\dotsc,z_6$ and the corresponding linear substitutions of the $x$ and $u$). Then, to find a covariant point $x_1,x_2,x_3$ for one of the permutations of the $z_1,\dotsc,z_6$, one only has to determine one more root of an easy cubic equation, namely to go to a fixed point of the connection. \\

In particular, Gordan succeeds in setting up a linear form of the desired kind, which is of degree 6 in the $z_1,\dotsc,z_6$. In the meantime, Mr. Coble showed by a systematic process that one need only go to the fourth degree \cite{Coble1911}. He sets up the associated cubic equation and then further sketches the course of the required algebraic calculation to determine the $z_1,\dotsc,z_6$. K. \\

I will conclude by mentioning the explanations given on p.491, footnote 10, referring to Kronecker's theorem. \\

First of all, for the sake of completeness, a few hints about my original proof of January, 1877. At that time, I had operated on the fact that all icosahedral forms, and also the tetrahedral form, have a direct degree. The circumstance is, of course, in turn, a consequence of the doubling of the number of homogeneous substitutions which I have emphasized in Abh. LIV, which was the actual reason for the proof. \\

Incidentally, I must go into more detail about the reference between Kronecker and myself from Easter 1881 [see p.10]. At that time, I asked Kronecker for a manuscript from 1861, from which I could copy the part that was suitable for me (the copy bears the date of March 23). In the proof of his theorem, Kronecker uses it exactly as I did later by anticipating L\"uroth's theorem on rational curves \cite{Luroth} and, from there on, the task is to form a rational function $\frac{\phi}{\psi}$ out of five free variables $x_0,\dotsc,x_4$, which is a linear transformation in the 60 alternating permutations of the $x_0,\dotsc,x_4$. He then encounters a strange lapse. Since Kronecker had not yet familiarized himself with the general notion of a group of linear substitutions (of one variable), he erroneously concludes that the 60 linear transformations in question ought to arise from the repetition of the same linear substitution; that is, from a cyclic equation of 60th degree, which (according to Galois theory) is of course impossible. At that time, Kronecker went quiet [on the subject]. \\

\newpage
In the lectures of 1885-86, this mistake is corrected. It is concluded that in the case of the alternating permutations of the free variables $x_0,\dotsc,x_4$, the polynomials $\phi$ and $\psi$ would have to be substituted in a linearly binary manner, and further, that such a binary behavior is already impossible if one of the permutations fixes an $x_i$ and cyclically permutes the others $x_j$'s. Even without the evidence of this impossibility, I still find an unnecessary complication. I showed above( p. 485) that even in the Klein four group, the impossibility in question arises. In order to come to a contradiction, Kronecker instead combines an operation of the Klein four group with the cyclic permutation of the $x_1,x_2,x_3$ - this is less transparent. \\

Aside from these secondary points, a complete consensus exists. There remains only a subjective difference, which I already discussed in detail on pages p.158-159 in the book on the icosahedron, but which I do not want to leave untouched here because of its importance. For the first time in his investigations into the solution of equations of the fifth degree, Kronecker begged to have a clear distinction between the natural irrationalities (which are rational functions of $x_0,\dotsc,x_4$) and the other irrationalities (which I call accessory). Incidentally, in his first communication of 1858 \cite{Brioschi1858LetterToHermite}
, he himself makes an unobjectionable use of an accessory square root. Is is only in the later work of 1861 \cite{Kronecker1861} that he believes that he should forbid the use of accessory irrationalities in the theory of equations altogether. In his 1885-86 lectures, he maintains this verdict:
\begin{quote}
    ...the use of accessory irrationalities is "algebraically worthless," because it "tears apart" the type.
\end{quote} 

\noindent 
In order to emphasize this demand, he calls it the "Abelian postulate." In contrast to other authors of similar thinking, I have explored as far as possible in my papers printed here above, as in the book on the icosahedron, the efficacy of using naturally occurring accessory irrationalities. \\

There is a principled difference in this thinking. I do not want to further emphasize that Abel continues to use the roots of unity $\epsilon$ in his investigations into the solution of the equations by radicals [which in the context of his considerations are also accessory irrationalities (see above, p. 486 footnote 8)], which incidentally Kronecker continues to do himself, because otherwise he would not be able to act on the connection of the equation of the fifth degree with the Jacobian equations of the sixth degree. Nor do I want to argue that it [the use of accessory irrationalities] is as  advantageous in the theory of numbers as it is in function theory, because of the simplicity of the higher-order algebraic relations in transcendental fields. I only want to emphasize the fundamentals. 

\np
When presented with new phenomena (such as the efficacy of the accessory irrationalities), should we we stop developing these ideas to align with our current conceptions, or rather push back against our narrow, systematic ways of thinking and pursue the new ideas in an unbiased way? Should one be a dogmatist or a natural scientist, and endeavor to keep learning from new ideas?\\

There is nothing special to be inferred from Kronecker's original notes, which Mr. Hensel sent to me. These are mainly 23 unilaterally described folios, of which 1-10 refer to the work of 1858 and 11-23 to those of 1861. It is worth noting that the passages I copied in 1881 are missing. On the back of the pages 17-18, there are bills with fifth roots of unity, by virtue of which Kronecker has evidently been convinced that the $G_{60}$ of broken icosahedral substitutions really do exist. \\

The criticism which I then apply to the transmitted material is intended to reflect the high position which I have given to Kronecker's investigations on the equations of the fifth degree in the above reprinted essays, especially in the historical account of the book on the icosahedron (\cite{KleinLectures}, see p.141-161). Kronecker first found the path which leads into the fundamental questions of the theory, only he did not finish it at first and later, at least formally, refused to accompany others on the way forward. K.

\np
\section{Appendix A: Formal Bibliography}

\bibliography{bibliography}
\bibliographystyle{plain}

\end{document}